\documentclass[draft]{elsarticle}

\usepackage{amssymb}
\usepackage{amsmath}
\newcommand{\R}{\mathbb R}

\newcommand \del     \partial

 \newcommand{\C}{\mathbb C}

\newcommand{\clg}[1]{{\mathcal{#1}}}

\newcommand{\ol}[1]{{\overline{#1}}}

\newcommand \AN      {\clg{AN}}

\def\XXint#1#2#3{{\setbox0=\hbox{$#1{#2#3}{\int}$}
\vcenter{\hbox{$#2#3$}}\kern-.5\wd0}}

\newtheorem{theorem}{Theorem}[section]
\newtheorem{proposition}[theorem]{Proposition}
\newtheorem{remark}[theorem]{Remark}
\newtheorem{lemma}[theorem]{Lemma}
\newtheorem{corollary}[theorem]{Corollary}
\newtheorem{definition}[theorem]{Definition}
\newtheorem{example}[theorem]{Example}
\newproof{pf}{Proof}
\journal{}
\begin{document}

\begin{frontmatter}


\title{Operators that achieve the norm}
\author[xc]{X.~Carvajal\fnref{fn1}}
\ead{carvajal@im.ufrj.br} 
\author[xc]{W.~Neves\fnref{fn2}}
\ead{wladimir@im.ufrj.br} 
\fntext[fn1]{Partially supported by FAPERJ, Brazil, through the
grant E-26/111.564/2008 entitled {\sl "Analysis, Geometry and
Applications"}, by Pronex-FAPERJ through the grant E-26/ 110.560/2010 entitled \textsl{%
"Nonlinear Partial Differential Equations"}, and by the National
Council of Technological and Scientific Development (CNPq),
Brazil, by the grant 303849/2008-8.} 
\fntext[fn2]{Partially supported by FAPERJ, Brazil, through the
grant E-26/111.564/2008 entitled {\sl "Analysis, Geometry and
Applications"}, and by
Pronex-FAPERJ through the grant E-26/ 110.560/2010 entitled \textsl{%
"Nonlinear Partial Differential Equations"}.
} 

\address[xc]{Institute of Mathematics, Federal University of Rio de Janeiro, C.P. 68530, Cidade Universit\'aria, 21945-970, Rio de Janeiro, Brazil}
\begin{abstract}
In this paper we study the theory of operators on complex Hilbert
spaces, which achieve the norm in the unit sphere. We prove
important results concerning the characterization of the $\AN$
operators, see Definition \ref{NDD}. The class of $\clg{AN}$
operators contains the algebra of the compact ones.
\end{abstract}
%
\begin{keyword}
Hilbert spaces, bounded operators, spectral representation.
%
\end{keyword}
\end{frontmatter}

\section{Introduction} \label{IN}

The theory of invariant subspaces of operators remains nowadays an
activity research area in functional analysis. It is still an open
question whenever an arbitrary operator on an infinite Hilbert
space has an invariant subspace, unless the trivial ones, that is
to say the whole space and the zero one.

We shall be concentrated on this article on a class of bounded
linear operators on complex Hilbert spaces, or on a subspace of
it, which attains his norm on the unit sphere. Here by a subspace,
we are always saying a closed subspace, and it is called invariant
under an operator, when such an operator maps the subspace into
itself.

In fact, the investigation of invariant subspaces for an operator
is the first step to understand better the structure of the
operator. For instance, the structure theorems on finite
dimensional case, the Jordan Decomposition Theorem and the
Spectral Theorem for normal operators illustrate, in particular,
decompositions of the Hilbert space in invariant subspaces.
Although, even for the finite dimensional set, the question of
finding a complete set of invariant subspaces for any operator on
a Hilbert space is arduous, unless for the self-adjoint and the
normal ones. In fact, the set of unitary invariants for a normal
matrix is obtained by its spectrum, that is counting the
multiplicities.

On the other hand, no general theory exists for a general operator
in the infinite Hilbert space set. For the infinite dimensional
case, the problem of how to count the spectral multiplicities for
self-adjoint operators was first introduced by Hellinger
\cite{EH}. Moreover, the question of counting the spectral
multiplicities for normal operators on infinite Hilbert spaces
could be studied in some different ways. For instance, in the
context of $C^*$-algebras.

\subsection{Purpose and some results} \label{PR}

In this paper we describe some new results on bounded operators in
Hilbert spaces, which achieve the norm. We study the operators
that satisfy the $\clg{N}$ and $\clg{AN}$ properties, defined
respectively in Definition \ref{ND} and Definition \ref{NDD}. This
class of operators contains, for instance, the compact ones.

Let $H$, $J$ be complex Hilbert spaces and $\clg{L}(H,J)$ the
Banach space of linear bounded operators from $H$ to $J$. We
emphasize the case that will appear most frequently later, namely
$\clg{L}(H,H)= \clg{L}(H)$.
Further, we recall that, the space $\clg{L}(H,J)$ is a Banach
space with the norm
\begin{equation} \label{ON}
   \|T\|= \sup_{\|x\|_{H} \leq 1}\|T x\|_{J}
\end{equation}
and, it is well known that, if $H$ has finite dimension, then the
closed unit ball in $H$ is compact (Heine-Borel Theorem) and the
above ``supreme'' is a maximum. In other words, if the dimension
$H$ is finite and $T \in \clg{L}(H,J)$, then there exists an $x$
in the closed unit ball in $H$ (indeed in the boundary, i.e. the
unit sphere), such that
$
   \|T\|= \|T x\|_{J}.
$
Although, in infinite dimension spaces this important property is
lost, albeit it remains true for the compact operators. Indeed, by
definition \eqref{ON} there exists a sequence $(x_n)_{n=1}^\infty$
in the unit closed ball of $H$, such that
$$
  \lim_{n \to \infty} \|T \, x_n \|_J = \|T\|.
$$
According to Banach-Alauglu Theorem,  an element $x$ exists in the
unit closed ball, and a subsequence $(x_{n_k})_{k=1}^\infty$ of
$(x_n)_{n=1}^\infty$, such that, $x_{n_k}$ converges weakly to
$x$. Therefore, since $T$ is a compact operator
$$
  \|T \, x_{n_k}\|_J \to \|Tx\|_J.
$$
It follows that, $\|T\|= \|T x\|_{J}$. Hence any compact operator
achieves the norm on the unit sphere. Then it is natural
to propose the following question: How to characterize the
operators which achieve the norm on the unit sphere?
Trying to answer it,
and also to study this interesting problem, we begin with the
following

\begin{definition}
\label{ND} An operator $T \in \clg{L}(H,J)$ is called to satisfy
the property $\clg{N}$, when there exists an element $x$ in the
unit sphere, such that
$$
  \|T\|= \|T x\|_{J}.
$$
\end{definition}

On the other hand, the restriction of a compact operator to a
subspace is a compact operator. Consequently, one can easily
observe that if $(x_n)_{n=1}^\infty$ is a sequence in the
intersection of the unit closed ball and a subspace $M$ of $H$,
and converges weakly, say $x_n \rightharpoonup x$, then $x \in M$,
since $M$ is closed.
Therefore a compact operator satisfies the following
generalization of the property $\clg{N}$. That is to say that a
compact operator achieves absolutely the norm (Absolutely Norm).
This suggests the following

\begin{definition}
\label{NDD} We say that $T \in \clg{L}(H,J)$ is an $\AN$ operator,
or to satisfy the property $\AN$, when
for all subspace $M \subset H$ $(M \neq \{0\})$, $T|_M$ satisfies
the property $\clg{N}$.
\end{definition}

We recall that by a subspace, we always mean a closed subspace,
hence on the definition quoted above $M$ is always closed.
Furthermore, it is quite normal to ask ourselves
 the following question: How to characterize the $\AN$ operators,
i.e. the ones that achieve absolutely the norm on the unit sphere?

One address the work of Bernard Chevreau, see \cite{BC}, who was
the first, as it is the knowledge of the authors, to introduce
some of these questions, when he was developing the canonical
writing form of compact operators without the use of spectral
properties, see \cite{IC}.

\subsection{Notation and background} \label{NB}

At this point we fix the functional notation used in this paper,
and recall some well known results from function analysis, see
\cite{FHHSPZ}, \cite{WR}.

By $(H,\langle .,. \rangle)$ we always denote a complex Hilbert
space, $S$ will denote the unit sphere in $H$ and $B$ the closed
unit ball in $H$.

If $T \in \clg{L}(H,J)$, the adjoint operator of $T$ is denoted by
$T^* \in \clg{L}(J,H)$, which satisfies $\|T^*\|=\|T\|$. An
operator $P \in \clg{L}(H)$ is called positive, when $\langle P \,
x,x \rangle \geq 0$, for all $x \in H$.
Given an operator $T \in \clg{L}(H,J)$, we denote by $P_T$, the
unique operator called the positive square root of $T^* T$, that
is, $\langle P_T \, x,x \rangle \geq 0$ for all $x \in H$ and
$P_T^2= T^* T$. Moreover, for $T \in \clg{L}(H)$ we recall the
polar decomposition of $T$, that is $T= U P$, where $U$ is a
unitary operator $(U^*= U^{-1})$ and $P \geq 0$.

As usual, if $x,y \in H$, then $x \perp y$ means that $x$ is
orthogonal to $y$, i.e. $\langle x,y \rangle= 0$. Additionally, if
$M \subset H$, we define
$$
  M^\perp:= \{ x \in H : \langle x,y \rangle= 0, \text{for
  all} \; y \in M \},
$$
that is the orthogonal complement of $M$, which is a (closed)
subspace of $H$. If $M$ is a subspace of $H$, hence closed by
assumption, then we could write
$
  H= M \oplus M^\perp.
$ For $x \in H$, we denote by $\C \,x$ the one-dimensional
subspace spanned by $x$, and by $x^\perp$ the orthogonal
complement of it.

\section{Property $\clg{N}$} \label{PN}


\begin{definition}\label{d1}
An operator $U \in \clg{L}(H,J)$ is called a partial isometry if
there exists a closed subspace $M$ of $H$, such that
$$
\forall x \in M,\,\,\, \|Ux\|=\|x\| \quad \textrm{and} \quad
\forall x \in M^{\bot}, \,\, Ux=0.
$$
Moreover, the orthogonal complement of the kernel of U is called
the initial domain, and the range of U the final domain.
\end{definition}

\smallskip
\begin{remark}
\label{PI} Let $U$ be a partial isometry with initial domain $M$
and final domain $M'$. Then, we have:

\medskip
{\bf1.} $U U^*$ is equal to the orthogonal projection on $M$,
$($denoted $P_M)$.

\medskip
{\bf 2.} $U^*$ is a partial isometry with initial domain $M'$.
\end{remark}

It is clear from Definition \ref{d1}, that if $U$ is a partial
isometry then $U$ satisfies the property $\clg{N}$. Additionally,
we observe that a linear combination of compact operators is also
a compact operator and therefore satisfies $\clg{N}$. Although,
the next example shows that this does not happen with the partial
isometry. Therefore, the set of operators that achieve the norm
does not form an algebra.

\begin{example} Let $\{e_j\}$ be an orthonormal base in $l^2$ and
$a \in (0,1]$.
Let $(a_j)$, $(b_j)$ be two sequences of real numbers, such that
$$
 0<a_1<a_2< \cdots <a, \quad a_j \nearrow a, \quad a_j^2+b_j^2=1.
$$
Now, let $T$ be the unitary operator given by
$
  Te_j:= \lambda_j \; e_j,
$
where $\lambda_j=a_j + i b_j$, $(j=1,2 \ldots)$. It is not
difficult to see that $(T+I)$ does not satisfy $\clg{N}$. In fact,
we have for each $x \in l^2$,
\begin{align*}
\left\|(T+I)x\right\|^2=&\sum_j\left\{|\lambda_j|^2+\lambda_j+\overline{\lambda_j}+1\right\}|x_j|^2
=\sum_j 2(1+a_j)|x_j|^2\\
< & \sum_j 2(1+a)|x_j|^2= 2(1+a)\left\|x\right\|^2.
\end{align*}
Moreover, $\left\|(T+I) \, e_j\right\|=\sqrt{2(1+a_j)}$, hence
$$
  \begin{aligned}
  \|(T+I)\|& \geq \sup_{j} \left\|(T+I) \,e_j\right\|
  = \lim_{j \to \infty} \sqrt{2(1+a_j)}= \sqrt{2(1+a)}.
  \end{aligned}
$$
Consequently, we have for any $x \in S$ that
$$
\left\|(T+I)\right\|= \sqrt{2(1+a)} \quad \textrm{and}
\,\,\left\|(T+I)x\right\|<\left\|(T+I)\right\|.
$$
\end{example}
\begin{proposition}\label{l1}
Let $T \in \clg{L}(H)$ be a self-adjoint operator. Then, $T$
satisfies $\clg{N}$ if, and only if $\|T\|$ or $-\|T\|$ is an
eigenvalue of $T$.
\end{proposition}
\begin{pf}1. The case when $\pm \|T\|$ is an eigenvalue of $T$ is
obvious. Indeed, if $x \in H$ is an eigenvector associated to $\pm
\, \|T\|$, then
$
  \|T x_0\|= \|T\|,
$
where $x_0:= x/\|x\|$.

2. Now, assume that there exists an element $x_0 \in S$, such that
$\|Tx_0\|=\|T\|$. Furthermore, without lost of generality, we can
suppose that $\|T\|=1$. We have
$$
\langle (I-T^2)x_0,x_0\rangle=\|x_0\|^2-\|Tx_0\|^2=0,
$$
and since $(I-T^2) \ge 0$ is a positive operator, it follows that
$(I-T^2)x_0=0$. Consequently, we must have
$$
  (I+T)(x_0 - T x_0)=0 \quad \text{or} \quad
  (I-T)(x_0 + T x_0)= 0.
$$
Let us look at the former case, the second one is analogue. If
$(x_0 - T x_0)= 0$, then we are done, since $\|T\|=1$ is an
eigenvalue of $T$ with the corresponded unitary eigenvector $x_0$.
On the other hand, if $z= (x_0 - T x_0) \neq 0$, let $z_0:= z/
\|z\|$. It follows that, $T z_0= -z_0$, that is, $z_0$ is a
unitary eigenvector associated to eigenvalue $-\|T\|=-1$.
\end{pf}

According to the proposition quoted above if $P \in \clg{L}(H)$ is
a positive operator and there exists an element $x_0 \in S$, such
that $\|Px_0\|=\|P\|$, then
\begin{equation}
\label{PX0} P x_0= \|P\| x_0.
\end{equation}
Likewise,
since $T$ satisfies $\clg{N}$ if, and only if $P_T$ satisfies
$\clg{N}$. Indeed,
\begin{equation}\label{rq}
\|T\|=\| P_T \| \quad \textrm{and}\quad \forall x \in H, \quad
\|Tx\|=\| P_T x\|,
\end{equation}
hence we have the following

\begin{corollary}\label{p1}
%
An operator $T \in \clg{L}(H,J)$ satisfies $\clg{N}$ if, and only
if $\|T\|$ is an eigenvalue of $P_T$.
\end{corollary}

Hereupon, we demonstrate the relation of the $\clg{N}$ condition
and the adjoint operator.

\begin{proposition}
\label{NCTAT}
 Let $T \in \clg{L}(H,J)$, then $T$ satisfies the condition $\clg{N}$ if,
 and only if the adjoint operator $T^*$ satisfies $\clg{N}$.
\end{proposition}

\begin{pf}
Since we are considering the bounded case, we have $T^{**}=T$.
Therefore, it is sufficient to prove one direction. Assume that
$T$ satisfies $\clg{N}$ condition. Then, $\|T x_0\|= \|T\|$ with
$x_0 \in S$. Now, set $P_T$ the positive square root of $T^*T$.
Hence by \eqref{PX0} and Corollary \ref{p1}, $P_T x_0=\|T\| x_0$,
and thus $T^*Tx_0=\|T\|^2 x_0$. Consequently, we have
$$
\left\|T^*\left(T\frac{1}{\|T\|}x_0\right)\right\|=\|T\|,
$$
and $T\left(\frac{1}{\|T\|}x_0\right) \in S$, which means that
$T^*$ carries $\clg{N}$ out.
\end{pf}

\begin{proposition}\label{inv} If $T \in \clg{L}(H,J)$ satisfies $\clg{N}$,
that is, if there exists an element $x_0 \in S$, such that $\|T
x_0\|=\|T\| $, then
\begin{equation}\label{ort}
T(x_0^\bot) \subset \big(Tx_0\big)^\bot.
\end{equation}
\end{proposition}

\begin{pf} By \eqref{rq}, we have $ \|P_T x_0\|=\|P_T\| $ and by
\eqref{PX0}, \eqref{rq} and Corollary \ref{p1}, we obtain $P_T \,
x_0=\|P_T\| x_0$. Therefore, if $y$ is orthogonal to $x_0$, then
$$
  T^*Tx_0=P_T^2 \, x_0= \|T\|^2 x_0
$$
and $y$ are orthogonal. Hence, $Ty$ is orthogonal to $Tx_0$.
\end{pf}

\begin{remark} \label{RRT}
If $T \in \clg{L}(H)$ is a positive operator and $x_0 \in S$ is
such that, $\|T x_0\|=\|T\|$, then
$\C x_0$ reduces $T$.
\end{remark}

It is not clear that all positive operators carry $\clg{N}$ out,
specially the self-adjoint ones. Then, we close this section with
an example of a positive operator, which does not satisfy the
$\clg{N}$ condition.
\begin{example}\label{ex2}
Let $T \in \clg{L}(l^2)$ be a positive operator, defined by
$$
  T: l^2 \to l^2, \quad (x_j) \mapsto (\lambda_j x_j),
$$
where $0<\lambda_1<\lambda_2< \ldots < \lambda$, with $\lambda_j
\nearrow \lambda$ and $\lambda < \infty$.

We have $T \ge 0$, $\|T\|=\lambda$, but for each $x \in l^2$ as
$Tx \neq \lambda x$, we conclude that $\|T\|$ is not an
eigenvalue. Consequently, by Proposition \ref{l1}, $T$ does not
satisfy the $\clg{N}$ condition.
\end{example}

\subsection{The numerical range relation}

\begin{definition}
Let $T \in \clg{L}(H)$. The numerical range of $T$ is defined as
$$
  W(T):= \{\langle Tx,x \rangle \in \C; \, x \in S \}.
$$
\end{definition}

\begin{remark}
Concerning the Toeplitz-Hausdorff's Theorem, $W(T)$ is a convex
set.
\end{remark}

Now, if $T \in \clg{L}(H)$ is a self-adjoint operator, then by
straightforward calculation $\|T\|= \sup_{x \in S} |\langle Tx,x
\rangle|$.
%
%
Therefore, for $P \geq 0$, it follows that
\begin{equation}\label{normap}
  \|P\|= \sup_{x \in S} \langle
  Px,x \rangle=\sup W(P).
\end{equation}

\begin{definition}
Let $A \subset \C$ be a convex non-empty set.
A number $\alpha \in A$ is said to be an extreme point of $A$,
when $\alpha= t \, u+(1-t) \, v$, with $u,v \in A$ and $0<t<1$
implies, $\alpha=u=v$.
\end{definition}

\begin{lemma}\label{l2}
Let $T \in \clg{L}(H)$ be an self-adjoint operator. Then, $T$
satisfies $\clg{N}$ if, and only if $\|T\|$ or $-\|T\|$ is an
extreme point of the numerical range $W(T)$.
\end{lemma}
\begin{pf}
First, we consider the operator $P^{\pm}: =\|T\| I \pm T$, which
is positive. Hence for all $x \in S$, we get
\begin{equation}\label{eq1}
\langle Tx,x \rangle \le \langle Tx,x \rangle \pm \langle
P^{\pm}x,x \rangle= \pm \|T\|.
\end{equation}
If $T$ satisfies $\clg{N}$, then  by Proposition \ref{l1}, we
obtain that $\|T\|$ or $-\|T\|$ is an extreme point of the
numerical range $W(T)$, since that
$$
\langle Tx,x \rangle \le \|T\|, \quad \langle Tx,x \rangle \ge
-\|T\|, \quad \forall x \in S.
$$

Now, if $\pm \|T\|$ is an extreme point of the numerical range
$W(T)$, then there exists $x_0 \in S$, such that $\pm \|T\|=
\langle T x_0, x_0 \rangle$. The equality in \eqref{eq1} implies
that $\langle P^{\pm} x_0, x_0 \rangle=0$. Moreover, since the
operator $P^{\pm}$ is positive, $P^{\pm} x_0=0$. Therefore, we
conclude that
$
  T \, x_0 =\pm \|T\| \, x_0.
$
\end{pf}
%

\section{The $\clg{AN}$ operators} \label{ANP}

As already seen at the Section \ref{PR} any compact operator $T$
in $\clg{L}(H,J)$ is an $\AN$ operator. Indeed, if $M$ is any
closed subspace of $H$, then $T|_M$ is compact and therefore
satisfies $\clg{N}$. Consequently, the algebra of compact
operators carries $\clg{AN}$ out. Since an orthogonal projection
is a partial isometry, it follows that any projection satisfies
the properties $\clg{N}$. Although, it is not necessarily true
that each projection is an $\AN$ operator. In fact, let us study
the following

\begin{example}
\label{ENAN} Let $X$ be the subspace of $l^2$, such the elements
$x$ have the form
$$
  x= (x_1, x_2, x_2, x_3, x_4, x_4, x_5, \ldots)
$$
and $P$ is the projection on $X$, i.e., $P:l^2 \to l^2$,
$$
  P(x_1,x_2,x_3,\ldots)= \Big(x_1, \frac{x_2+x_3}{2},
  \frac{x_2+x_3}{2}, x_4, \frac{x_5+x_6}{2}, \frac{x_5+x_6}{2},
  x_7, \ldots \Big).
$$
Now, let $M$ be a subspace of $l^2$, defined as
$$
  M:= \{x \in l^2 : x= (x_1,x_1,x_2,x_2,x_2,x_3,x_3,x_3,x_4,x_4,x_4,\ldots)
  \}.
$$
It follows that, $M \cap X= \{0\}$. Set $P|_M \equiv T: M \to
l^2$, hence
$$
  T(x_1,x_1,x_2,x_2,x_2,\ldots)= \Big(x_1, \frac{x_1+x_2}{2},
  \frac{x_1+x_2}{2}, x_2, \frac{x_2+x_3}{2}, \frac{x_2+x_3}{2},
  x_3, \ldots \Big).
$$
For each $x \in M \cap S$, we compute the norm of $T x$. First, we
have
\begin{equation} \label{NX}
  1= \|x\|^2 = 2 x_1^2 + 3 \, \sum_{j=2}^\infty x_j^2.
\end{equation}
Then, it follows that
\begin{equation} \label{TNX}
  \begin{aligned}
  \|Tx\|^2&= \sum_{j=1}^\infty x_j^2 + 2 \, \sum_{j=1}^\infty \big(\frac{x_j+x_{j+1}}{2}\big)^2
\\
  &= x_1^2 + \sum_{j=2}^\infty x_j^2 + \frac{x_1^2}{2} + \sum_{j=2}^\infty
  \frac{x_j^2}{2} + \sum_{i=2}^\infty x_i^2 + \sum_{j=1}^\infty
  x_j \, x_{j+1}
  \\
  &= \frac{2}{3} + \frac{x_1^2}{6} + \sum_{j=1}^\infty
  x_j \, x_{j+1},
  \end{aligned}
\end{equation}
where we have used \eqref{NX}. Let $\{s^n\}_{n=1}^\infty$ be a
sequence contained in $M \cap S$, i.e. for each $n=1,2,\ldots$,
$s^n \in M \cap S$,
$$
  s^n=
  (s^n_1,s^n_1,s^n_2,s^n_2,s^n_2, \ldots, s^n_n,s^n_n,s^n_n, \ldots),
  \quad \|s^n\|=1,
$$
$s^n$ defined as
\begin{equation} \label{SNJ} s^n_j= \left \{
\begin{aligned}
      \frac{1}{\sqrt{3 (n-1) + 2}} \quad j&= 1,\ldots,n,
      \\
      0 \quad \quad \quad \quad j&>n.
\end{aligned}
\right.
\end{equation}
We claim that, $\|T s^n\| \to 1$ as $n \to \infty$. Indeed, from
\eqref{TNX}
$$
    \begin{aligned}
  \|T s^n\|^2&= \frac{2}{3} + \frac{(s^n_1)^2}{6} + \sum_{j=1}^\infty
  s^n_j \, s^n_{j+1}
  = \frac{2}{3} + \frac{6n - 5}{6(3(n-1)+2)}.
  \end{aligned}
$$
Then, we obtain
$$
  \lim_{n \to \infty} \|T s^n\|^2= \frac{2}{3} + \frac{6}{18}= 1.
$$
Consequently, $\|T\|= \|P\|= 1$.

Now, let us show that $T$ does not satisfy $\clg{N}$. If there
exists $x \in S \cap M$, such that
$
  1= \|Tx\|= \|Px\|,
$ then $1$ is an eigenvalue of $P$ associated to the vector $x$.
Consequently, $Px= x$, which is a contradiction, since we have
assumed $M \cap P(l^2)= \{0\}$ and we have $\|x\|=1$.
\end{example}

Therefore, we have proved the following result.

\begin{lemma}
\label{CE1} Let $P$ be an orthogonal projection. Then $P$ does not
necessarily satisfy $\clg{AN}$ property.
\end{lemma}

The next proposition will be used as a proof of the next theorem,
but it is important by itself.

\begin{proposition}
\label{PCOM} Let $R$ be an isometry on $H$ and $T \in \clg{L}(H)$
an $\clg{AN}$ operator. Then, $TR$ and $RT$ satisfy the property
$\clg{AN}$.
\end{proposition}

\begin{pf}
Let $M$ be a subspace of $H$. First, clearly the composition $RT$
satisfies $\clg{AN}$. Indeed, if $x_0 \in M \cap S$ is such that,
$\|T\|= \|T x_0 \|$, then
$$
  \|R \, T \|= \|T\|= \|T x_0\|= \|R \, T x_0\|.
$$
Now, let us show that $TR$ also satisfies $\clg{AN}$. Since $M$ is
a closed subspace of $H$ and $R$ an isometry, hence $R(M)$ is also
a closed subspace of $H$. Moreover, we have
$$
  \| T \, R|_M \|= \| T |_{R(M)} \|.
$$
Since $R(M)$ is closed and $T$ satisfies $\clg{AN}$, there exists
$Rz \in R(M)$, with $\|Rz\|= \|z\|= 1$, such that
$
  \|T \, R|_M \|= \|T(Rz)\|.
$
\end{pf}

Subsequently, we recall a well known definition for equivalent
operators.

\begin{definition}
\label{EQO} The operators $T \in \clg{L}(H)$ and $S \in
\clg{L}(J)$ are called unitary equivalents, when there exists a
unitary operator $U$ on $\clg{L}(J,H)$, such that
$$
  U^* \, T \, U= S.
$$
\end{definition}
In fact, if $T$ and $S$ are unitary equivalents, then there is no
criterion based only on the geometry of the Hilbert space, in such
a way that, $T$ could be distinguished from $S$. Therefore, since
$T$ and $S$ are abstractly the same operator, it is natural to
conjecture that some characteristic endowed by $T$ must be
satisfied by $S$, and vice-versa.

\begin{theorem}
\label{UEO} Let $T,S$ be two unitary equivalent operators. Then,
$T$ is $\clg{AN}$ operator if, and only if $S$ is an $\clg{AN}$
operator.
\end{theorem}

\begin{pf}
Assume that $U$ is a unitary operator such that $U^\star \, T \,
U= S$, hence $TU= US$. Since $U$ is an isometry, by Proposition
\ref{PCOM} if $T$ satisfies $\clg{AN}$, then $TU$ satisfies
$\clg{AN}$. Moreover, it follows that, $US$ also satisfies
$\clg{AN}$. Once more, conforming to Proposition \ref{PCOM}, we
have that $S$ satisfies property $\clg{AN}$.
\end{pf}

\begin{remark} \label{TPT}
Given $T \in \clg{L}(H,J)$, we recall that $P_T$ was defined as
the positive square root of $T^*T$. Therefore, $T$ satisfies
$\clg{AN}$ if, and only if $P_T$ satisfies $\clg{AN}$, see
\eqref{rq}. Consequently, it is enough to establish the condition
$\clg{AN}$ for positive operators.
\end{remark}

\begin{proposition}
An operator $T \in \clg{L}(H,J)$ satisfies the property $\clg{AN}$
if, and only if, for all orthogonal projection $Q \in \clg{L}(H)$,
$TQ$ satisfies $\clg{N}$.
\end{proposition}

\begin{pf}
Let $M$ be a closed subspace of $H$ and $Q$ an orthogonal
projection on $M$. Then, we have
$
  \|TQ\|= \|T|_M\|.
$
\end{pf}

\begin{lemma}
\label{IFR} Let $R \in \clg{L}(H)$ be an operator of finite rank.
Then $I + R$ is an $\clg{AN}$ operator.
\end{lemma}

\begin{pf}
Suppose that $\dim R(H)=n$. Since $R$ is finite rank, we could
write
$$
  R x = \sum_{j=1}^n \lambda_j \langle x, e_j\rangle \, e_j,
$$
where $\{e_j\}_{j=1}^n$ is an orthonormal set of $H$ and
$\lambda_j \ge 0$, $(j=1,2, \dots n)$. Let $M_n$ be the subspace
generated by $\{ e_1,\ldots,e_n \}$, thus
$
  H= M_n \oplus M_n^\perp.
$
Moreover, for any $x \in H$, $x= x_1 + x_2$, such that
$$
  x_1= \sum_{j=1}^n \langle x,e_j \rangle \, e_j
  \quad \text{and} \quad
  x_2= \sum_{\alpha \in A} \langle x,e_\alpha \rangle \, e_\alpha,
$$
where $\{e_\alpha\}_{\alpha \in A}$ is an orthonormal basis of
$M_n^\perp$, $\langle e_\alpha,e_j \rangle= 0$, for all
$j=1,\ldots,n$, $\alpha \in A$ and $\langle e_\alpha,e_\beta
\rangle= \delta_{\alpha \beta}$ for each $\alpha, \beta \in A$.
Now, define $T:= I + R$, then for each $x \in H$,
$$
  Tx= \sum_{j=1}^n \langle x,e_j \rangle \, e_j
  + \sum_{\alpha \in A} \langle x,e_\alpha \rangle \, e_\alpha
  + \sum_{j=1}^n \lambda_j \langle x, e_j\rangle \, e_j.
$$
Consequently, for each $x \in S$,
$$
  \|T x\|^2= 1 + \sum_{j=1}^{n}
  \big(\lambda_j^2 + 2 \, \lambda_j \big)
  |\langle x, e_j \rangle|^2.
$$
Therefore, if $P$ is the finite range projection on $M_n$, then $
  \left\|T P x\right\|=\left\|T x\right\|,
$ for any $x \in S$, and as $T P$ has finite range and therefore
satisfies $\clg{AN}$, then $T$ satisfies $\clg{AN}$.
\end{pf}

\begin{lemma}\label{LPQRN} Let $H$ be
a separable Hilbert space. If $P,Q \in \clg{L}(H)$ are two
orthogonal projections such that, the dimension of their ranks and
null spaces are infinite, then $P$ and $Q$ are unitary equivalent.
\end{lemma}

\begin{pf}
Since the rank and the null space of a projection are subspaces,
there exist unitary operators $U_1: P(H) \to Q(H)$ and $U_2:
\clg{N}(P) \to \clg{N}(Q)$. Now, we define $U: H \to H $, such
that
$$
  U|_{P(H)}=U_1 \quad \text{and} \quad U|_{\clg{N}(P)}=U_2.
$$
Hence it is clear that $U$ as defined above is a unitary operator.
Moreover, if $x \in H$, then $x=x_1+x_2$ where $x_1 \in P(H)$ and
$x_2 \in \clg{N}(P)$. From the definition of $U_1$ and $U_2$, we
have
$
  Q U x=U x_1=U P x.
$
Therefore, $P$ and $Q$ are unitary equivalents.
\end{pf}

\begin{theorem}
\label{PONR} Let $Q \in \clg{L}(H)$ be an orthogonal projection.
Then, $Q$ satisfies $\clg{AN}$ property
 if, and
only if the dimension of the null space or the dimension of the
rank of $Q$ is finite.
\end{theorem}

\begin{pf}
If $\dim  Q(H)< \infty$, then $Q$ is compact, hence it satisfies
$\clg{AN}$. Now, assume $\dim \clg{N}(Q)< \infty$. Then, we have
$Q= I - P$, where $P$ is a projection with finite rank. Therefore,
$Q$ satisfies $\clg{AN}$.

On the other hand, if $\dim Q(H)= \dim \clg{N}(Q)= \infty$, we
consider two cases:

i)  H separable. In this case, by Lemma \ref{LPQRN}, we have that
$Q$ is unitary equivalent to the orthogonal projection of Example
\ref{ENAN}, which does no satisfy $\clg{AN}$ condition.
Consequently, by Theorem \ref{UEO} $Q$ does not satisfy $\clg{AN}$
either.

ii) H is not separable. If $Q(H)$ is countable, we take $H_1$ be a
separable Hilbert space such that $Q(H) \subset H_1$ and $\dim
(Q(H)^{\bot} \cap H_1)=\infty $. Thus, we have that $Q|_{H_1}$ is
an orthogonal projection, $Q|_{H_1} \in \clg{L}(H_1)$ and by the
separable case (i), $Q|_{H_1}$ does not satisfy $\clg{AN}$.
Therefore, we have that $Q$ does not satisfy $\clg{AN}$ either.

Now, if $Q(H)$ is not countable, let $H_1 \subset Q(H)$ be a
countable subspace and $Q_1$ be an orthogonal projection on $H_1$.
Furthermore, let $N_1$ be an infinite countable subset
of $Q(H)^{\bot}$, then $H_2 =H_1 \oplus N_1$ is a separable
Hilbert space. Conforming with the separable case (i), it follows
that $Q_1|_{H_2}\in \clg{L}(H_2)$ is an orthogonal projection in
$H_2$ which does not satisfy $\clg{AN}$, since $\|Qx\|=\|Q_1x\|$
for all $x \in H_2$. Consequently, neither $Q$ satisfy $\clg{AN}$.
\end{pf}

\begin{definition}
\label{COISO} Let $B \in \clg{L}(H)$. The operator $B$ is called a
co-isometry, when $B^*$ is an isometry.
\end{definition}

As stated in Remark \ref{PI}, it follows that:

\begin{remark} One observer that, a co-isometry $V$
is a partial isometry with initial domain $V^*(H)$. Thus, if $V^*$
is not a unitary operator, then $\clg{N}(V) \neq \{0\}$.
\end{remark}

Now, we observe that if $B \in \clg{L}(H)$ is a co-isometry, then
$$
  BB^*= I, \quad \quad \clg{N}(B^*B)= \clg{N}(B).
$$
Furthermore, we have the following characterization for the
$\clg{AN}$ property.

\begin{proposition}
\label{PCOISO} Let $B$ be a co-isometry. Then, $B$ satisfies the
property $\clg{AN}$ if, and only if $\clg{N}(B)$ has finite
dimension.
\end{proposition}

\begin{pf}
Since $B^*$ is an isometry, by Proposition \ref{PCOM} $B$
satisfies $\clg{AN}$ condition if, and only if $B^* B$ satisfies
$\clg{AN}$. Now, we observe that $B^*B$ is an orthogonal
projection, hence by Theorem \ref{PONR}, it follows that $B$
carries $\clg{AN}$ out if, and only if $\dim \clg{N}(B^*B)<
\infty$ or $\dim B^*B(H)< \infty$. Additionally, we have
%
$$
  B^*B(H)= \clg{N}(B^*B)^\perp= \clg{N}(B)^\perp= B^*(H).
$$
Consequently, $B$ satisfies $\clg{AN}$ if, and only if $\dim
\clg{N}(B)< \infty$ or $\dim B^*(H)< \infty$. Likewise, since
$B^*$ is an injective operator, if $B^*(H)$ has finite dimension,
then $\clg{N}(B)$ also has finite dimension.
\end{pf}

Contrarily to Proposition \ref{NCTAT}, and according to the
proposition quoted above, it follows that:

\begin{remark}
Even if $V$ is an $\clg{AN}$ operator, it does not follow
necessarily that $V^*$ carries $\clg{AN}$ out.
\end{remark}

\begin{proposition}
\label{PSAN} Let $U \in \clg{L}(H)$ be a partial isometry with
initial domain $M$. Then, $U$ satisfies $\clg{AN}$ if, and only
if, the dimension of $M$ or the dimension of $M^\perp$ is finite.
\end{proposition}

\begin{pf}
For all $x \in H$, we have
$
  \|U^*U x\|= \|Ux\|.
$
Therefore, $U$ satisfies $\clg{AN}$ if, and only if $U^*U= P_M$
satisfies $\clg{AN}$. We can prove it applying Theorem \ref{PONR}.
\end{pf}

\begin{remark}
Let $P \in \clg{L}(H)$ be a positive operator. Then, $P$ is a
partial isometry with initial domain $M$ if, and only if $P$ is an
orthogonal projection on $M$.
%
Indeed, assume that $P$ is a partial isometry with initial domain
$M$. Hence $\|P\|= 1$ and for each $x \in M$, $\|P x\|= \|x\|$.
Since $P$ is a positive operator, it follows that
$$
  \forall x \in H, \quad \|P\|^2 \leq \|P\| \; \langle Px,x
  \rangle.
$$
Moreover, we have for all $x \in S$, $\langle Px,x
  \rangle \leq \|P\|= 1$.
Consequently, for each $x \in M \cap S$, $\langle Px,x \rangle=
1$.

Now, let $P$ be an orthogonal projection on $M$ and set $T:= I -
P$. Then, for each $x \in S \cap M$, we have
$
  \langle Tx,x \rangle= 0.
$
Since $T$ is a positive operator, we have for all $x \in M$, $Tx=
0$. Therefore, $T|_M= 0$, i.e. $P|_M= I$.
\end{remark}

\subsection{More characterization}

\begin{proposition}
Let $V,R \in \clg{L}(H,J)$ be respectively an isometry and a
finite rank operator. Then $V+R$ is an $\clg{AN}$ operator.
\end{proposition}
\begin{pf}
As $V^{*}V=I$, then $V+R=V+RV^{*}V=(I+RV^{*})V$. Since that
$RV^{*}$ is an operator of finite rank, by Proposition \ref{PCOM}
and Lemma \ref{IFR}, we conclude that $V+R$ satisfies $\clg{AN}$.
\end{pf}

\begin{proposition}\label{ANPVR}
Let $V,P \in \clg{L}(H)$ be respectively  a co-isometry and an
orthogonal projection, which are $\clg{AN}$ operators. If $R \in
\clg{L}(H)$ is an operator of finite rank, then
$
  V + R \quad \text{and} \quad P + R,
$
carry $\clg{AN}$ out.
\end{proposition}
\begin{pf}
Since $V$ is a co-isometry, similarly as above, we have
$V^{*}V=I$, thus $V+R=V(V^{*}R+I)$. Therefore, as $V^{*}R$ is an
operator of finite rank, by Proposition \ref{PCOM} and Lemma
\ref{IFR}, we conclude that $V+R$ is an $\clg{AN}$ operator.

Now, if $P$ satisfies $\clg{AN}$, then $P$ has finite rank, or we
could write $P=I-K$, where $K$ is a projection with finite rank.
Consequently, by Lemma \ref{IFR} $P+R$ is an $\clg{AN}$ operator.
\end{pf}

\begin{proposition}
\label{ANPIFR} Let $W \in \clg{L}(H,J)$ be a partial isometry with
initial domain $M$, which satisfies $\clg{AN}$ and $R \in
\clg{L}(H,J)$ an operator of finite rank. Then, $W+R$ is an
$\clg{AN}$ operator.
\end{proposition}

\begin{pf}
1. First, since $W$ carries $\clg{AN}$ out, hence, conforming to
Proposition \ref{PSAN}, we have $\dim M < \infty$ or $\dim M^\perp
< \infty$. We assume that $\dim M < \infty$, then $W$ is compact
and, consequently $W+R$ satisfies $\clg{AN}$ condition.

2. Now, we assume that $\dim M^\perp < \infty$. Set $n:= \dim
R(H)$, thus for any $x \in H$, we could write
$$
  R x= \sum_{j=1}^n \lambda_j \; \langle x, e_j \rangle \; \epsilon_j,
$$
where $\lambda_j>0$, $(j=1,\ldots,n)$ and $\{e_j\}_{j=1}^n$,
$\{\epsilon_j\}_{j=1}^n$ are orthonormal sets in $H$ and $J$
respectively. Let $m= \dim M^\perp$ and $\{\varphi_j\}_{j=1}^m$ be
an orthonormal basis for $M^\perp$. Moreover, for any $x \in H$,
it follows that
$
  \|W x\|^2= \|x\|^2 - \sum_{j=1}^m |\langle x,\varphi_j \rangle|^2.
$ Then, for any $x \in S$
$$
  \begin{aligned}
  \|(W+R)x&\|^2= \langle Wx+Rx,Wx+Rx \rangle
\\
  &= \|Wx\|^2 + \langle Wx,Rx \rangle + \langle Rx,Wx \rangle + \|Rx\|^2
\\
  &= 1 - \sum_{j=1}^m |\langle x,\varphi_j \rangle |^2 + 2\, \clg{R}e \sum_{j=1}^n \lambda_j
  \ol{\langle x,e_j \rangle} \, \langle x,W^{*}\epsilon_j \rangle
  + \sum_{j=1}^n \lambda_j^2  |\langle x,e_j \rangle |^2.
  \end{aligned}
$$
Now, let $N$ be the subspace generated by
$$
  \{e_j,W^{*}\epsilon_j,\varphi_k\}, \quad (j=1,\ldots,n, \,
  k=1,\ldots,m),
$$
and $Q$ the orthogonal projection on $N$. Hence, it follows that,
for all $x \in S$,
$$
  \| (W+R) Q x \|= \| (W+R) x \|.
$$
Therefore, as $(W+R) Q$ satisfies $\clg{AN}$, we conclude that
$W+R$ carry $\clg{AN}$ out.
\end{pf}

\begin{proposition}
Let $P_1,P_2 \in \clg{L}(H)$ be $\AN$ orthogonal projections. Then
$P_1 \pm P_2$, $P_1 P_2$ and $P_2 P_1$ satisfy the $\clg{AN}$
property.
\end{proposition}

\begin{pf}
In fact, the proof follows with the following remark. If $P$ is an
orthogonal projection, which satisfies $\clg{AN}$, then $P$ or
$I-P$ has finite rank. Therefore, if $P$ satisfies $\clg{AN}$ or
$P$ has finite rank, or we could write $P=I-K$, where $K$ is a
projection with finite rank.
\end{pf}

\begin{theorem}
Let $M$ be a subspace of $H$ and $T \in \clg{L}(H,J)$, such that
$$
  \|T|_M \| < \|T\|.
$$
Then, there exists an element $x \in S$, satisfying
$
  \|Tx\| = \|T|_M \|.
$
\end{theorem}

\begin{pf}
1. Since $\|T|_M \| < \|T\|$, given $\epsilon>0$, there exists $w
\in S$, such that $\|T|_M \| < \|T w\| + \epsilon$. Consequently,
we have
$$
  \|T|_M \|^2 \leq  \|T w \|^2= \langle T^* T w, w \rangle.
$$
Moreover, for some $y \in S$, $\|T y \| \leq \|T|_M \|$. Hence we
have
$$
  \langle T^* T y, y \rangle \leq \|T|_M \|^2 \leq \langle T^* T w, w \rangle.
$$

2. Now,
by the convexity of the numerical rank of the positive operator
$T^* T$, i.e. $W(T^*T)$,
there exists an element $x \in S$, satisfying
$
  \langle T^\star T x, x \rangle= \|T|_M \|^2.
$
Consequently, we obtain $\|T x\|= \|T|_M \|$.
\end{pf}

\begin{corollary}
Let $P \in \clg{L}(H)$ be a positive operator and $M$ a subspace
of $H$, such that $\|P|_M \| < \|P\|$. Then, there exists an
element $x \in S$, satisfying
$$
  \langle Px,x \rangle= \|P|_M\|.
$$
\end{corollary}
%
%

\begin{proposition}
Let $K \in \clg{L}(H)$ be a positive compact operator. Then, $K+I$
is an $\clg{AN}$ operator.
\end{proposition}

\begin{pf}
Since $K$ is compact and positive, it is also $K^2 + 2 \, K$.
Hence there exists $T \in \clg{L}(H)$ a compact and positive
operator, such that
$
  T^2= K^2 + 2 \, K.
$
Now, let $M$ be a subspace of $H$. Then, there exists $x \in S
\cap M =:S_M$, such that $\|T x\|^2= \| T|_M \|^2$, that is
$$
  \langle (K^2 + 2 \, K ) x, x \rangle= \sup_{z \in S_M} \| T z
  \|^2= \sup_{z \in S_M} \langle (K^2 + 2 \, K ) z, z \rangle.
$$
Consequently, we have
$$
  \begin{aligned}
   \| (K+I)|_M \|^{2} &= \sup_{z \in S_M}\| (K+I)z \|^{2}
   = 1 + \langle (K^2 + 2 \, K ) x, x \rangle
   = \| (K+I)x \|^{2}.
  \end{aligned}
$$
\end{pf}

On the other hand, we already know that some type of operator
satisfying $\clg{AN}$, such a sum with an operator of finite rank
already satisfy $\clg{AN}$. For instance, if $K \in \clg{L}(H)$ is
a positive compact operator, then $K+I$ satisfies $\clg{AN}$.
Therefore, observing the examples at the beginning of this
section, it is natural to ask, if $K+I+R$ satisfies $\clg{AN}$
condition, when $R \in \clg{L}(H)$ is an operator of finite rank.
In fact, this question was positively answered.

\begin{theorem}
\label{NCKIFR} Let $K \in \clg{L}(H)$ be a positive compact
operator and $R \in \clg{L}(H)$ an operator of finite rank. Then,
$K+I+R$ satisfies the $\clg{AN}$ condition.
\end{theorem}

\begin{pf}
Initially, we will prove that $K+I+R$ satisfies the $\clg{N}$
condition.

1. If $K$ has finite rank, then from Lemma \ref{IFR}, $K+I+R$
satisfies $\clg{N}$ condition.

2. Now, assume that $K$ does not have finite rank and set $T:=
I+K+R$. Furthermore, we can suppose $T$ positive and $R$
self-adjoint, since $T$ satisfies $\clg{N}$ if, and only if $T^*T$
does.
Indeed, we have $T^*T= I + \mathsf{K} + \mathsf{R}$, where
$\mathsf{K}= \big(K^*+K+ K^* K \big)$ is a positive compact
operator, and
$$
  \mathsf{R}= \big(R + R^* + R^* R + K^* R + R^* K
\big)
$$
is a self-adjoin finite rank operator. Therefore, we have for each
$x \in H$,
\begin{equation} \label{100}
  (K+R)x= \sum_{j=1}^\infty \lambda_j \, \langle x,e_j \rangle \,
  e_j,
\end{equation}
where $\{e_j\}$ is a orthonormal sequence, $|\lambda_j| \searrow
0$ and $\lambda_j \in \R$ for all $j \geq 1$. Consequently, for
each $x \in H$,
$$
  \langle (K+R)x,x \rangle = \sum_{j \geq 1} \lambda_j \; |\langle x, e_j
  \rangle|^2.
$$

3. We claim that, there exists a $k$, such that $\lambda_k > 0$.
Otherwise $0 \le K\leq -R$, but since $R$ has finite rank and $K$
is positive, it follows that $K$ has finite rank.
Indeed, let $n= \dim R(H)$. If $x \in R(H)^{\bot}$, then $ 0 \le
\langle Kx,x \rangle \le  -\langle Rx,x \rangle=0,$ therefore
$\langle Kx,x \rangle=0$ and since that $K \ge 0$, we have $Kx=0$
for all $x \in R(H)^{\bot}$. Moreover, if $R x= \sum_{j=1}^n
\gamma_j \, \langle x,\eta_j \rangle \,
  \eta_j$, then $K(H)$ is contained in the subspace generated by
  $\{K(\eta_1), K(\eta_2), \dots, K(\eta_n) \}$, which is a contradiction.

4. Let $E$ be the subspace generated by $\{e_j\}$, thus $H= E
\oplus E^\perp$. Consequently, for each $x \in H$, we have
\begin{equation} \label{Tx+}
\begin{aligned}
   T x= \sum_{j=1}^\infty (\lambda_j + 1) \; \langle x,e_j \rangle \,
  e_j + P x,
\end{aligned}
\end{equation}
where $P$ is an orthogonal projection on $E^\perp$.
Now, since $T$ is positive
$
  \langle Te_j,e_j \rangle =(\lambda_j + 1) \geq 0,
$
for all $j \geq 1$. From \eqref{Tx+}, it follows that
\begin{align*}
 \|T x\|^{2}=& \sum_{j=1}^\infty (\lambda_j + 1)^{2} \; |\langle x,e_j \rangle|^{2} \,
  +\| P x\|^{2}
  \\
  \le & \sup_{j} \{(\lambda_j + 1)^2 \} \sum_{j=1}^\infty  \; |\langle x,e_j \rangle|^{2} \,
  +\| P x\|^{2}
  = \sup_{j} \{(\lambda_j + 1)^2 \} \; \|x\|^{2},
 \end{align*}
where we have used that, $\sup_{j} \{\lambda_j + 1 \}^{2} \ge
(\lambda_k + 1)^{2}>1$. Then, we obtain
$$
  \| T \|= \sup_{j} \{\lambda_j + 1 \}.
$$
Now as $|\lambda_1| \ge \ldots \ge |\lambda_{k-1}| \ge\lambda_{k}
\ge |\lambda_{k+1}| \ge \ldots $, there exists an indices $i \leq
k$, such that
$
  \|T\|= \lambda_i + 1= \|T e_i\|,
$
which proves that $K+I+R$ satisfies the $\clg{N}$ condition.

Finally, let $M$ be a subspace of $H$, and define
$$
 T_1:= T|_M, \quad I_1:= I|_M, \quad K_1:= K|_M \quad
 \text{and} \quad R_1:= R|_M.
$$
Then, $T_1=I_1+K_1+R_1$. Since
$T_1$ satisfies the $\clg{N}$ condition if, and only if
$T_1^{*}T_1$ satisfies the $\clg{N}$ condition, similarly as
displayed above, we can suppose that $T_1 \in \clg{L}(M)$ is a
positive operator, $K_1 \in \clg{L}(M)$ is a positive compact
operator and $R_1 \in \clg{L}(M)$ is a self-adjoint finite rank
operator. Therefore, we can procedure as the initial part and
obtain that $T_1$ satisfies the $\clg{N}$ condition.
\end{pf}

\subsection{Structure of $\clg{AN}$ operators}

In this section, we analyze the structure of the positive
operators satisfying the $\clg{AN}$ condition.

\begin{example}
Let $H$ be a separable Hilbert space and $\{e_j\}$ a hilbertian
basis of it. Let $T \in \clg{L}(H)$ be defined as
$
  T x:= \sum_{j=1}^\infty \lambda_j \, \langle x,e_j \rangle \,
  e_j,
$
where $\lambda_j \searrow \lambda$. Then, $T$ is an $\clg{AN}$
operator. Indeed, set $V:= T - \lambda I$, therefore
$$
  V x= \sum_{j=1}^\infty (\lambda_j - \lambda) \, \langle x,e_j \rangle \,
  e_j.
$$
Since $\lambda_j \searrow \lambda$, $V$ is compact and positive,
hence $T= \lambda I + V$ is an $\clg{AN}$ operator.
\end{example}
Now
let $T \in \clg{L}(H)$ be defined as
\begin{equation}
\label{LOTD}
  Tx= \sum_{j=1}^{\infty}\lambda_j\langle x, \upsilon_j\rangle
\upsilon_j,
\end{equation}
where $\{\upsilon_j\}$ is an orthonormal set in $H$ (not
necessarily a basis) and $\lambda_j\searrow \lambda>0$.
%
If $\dim \clg{N}(T)< \infty$, then $T$ satisfies $\clg{AN}$
condition.
%
Indeed,
%
we write $T= K + \lambda P$, where for any $x \in H$
$$
  K x= \sum_{j=1}^{\infty} (\lambda_j - \lambda) \langle x, \upsilon_j\rangle
\upsilon_j, \quad \quad P x= \sum_{j=1}^{\infty} \langle x,
\upsilon_j\rangle \upsilon_j.
$$
Therefore, $K$ and $P$ are respectively a compact operator and
an orthogonal projection on $K(H) ~(=T(H))$. Now, since $\dim
\clg{N}(T)< \infty$, hence we have $P=I-R$, where $R$ is an
orthogonal projection of finite range on $\clg{N}(T)$. Thus
$$
  T=\lambda \Big(\frac{1}{\lambda} \, K + I - R \Big)
$$
satisfies $\clg{AN}$ by Proposition \ref{NCKIFR}.
%
%


\begin{definition}
Let $T_1,T_2 \in \clg{L}(H,J)$ be two operators, we say that $T_1$
and $T_2$ are mutually orthogonal, denoted $T_1 \bot T_2$, if for
any $x,y \in H $,
$
  \langle T_1 x,T_2 y\rangle=0.
$
Moreover, when $T_1 \bot T_2$ we write the sum of $T_1$, $T_2$ as
$T_1\oplus T_2$.
\end{definition}

\begin{theorem}\label{decomp+}
Let $T \in \clg{L}(H)$ be a positive operator, which satisfies the
$\clg{AN}$ condition .
Then $T$ has the following representation
$$
  T= \sum_{n \geq 1} \beta_n \; \upsilon_n
\otimes \upsilon_n
\oplus R_1,
$$
where $\{\upsilon_n\}_{n=1}^{\infty}$
is an orthonormal sequence of vectors in $H$, $\beta_n \searrow
\beta \ge 0$ and  $\beta \ge \|R_1\|$.
\end{theorem}

\begin{pf} It is enough to show the representation.
First, since $T \in \clg{L}(H)$ is a positive $\AN$ operator,
there exists an element $\upsilon_1 \in S$, such that $T \upsilon
_1= \|T\| \upsilon_1$. Let $H_1$ be the one-dimensional subspace
of $H$ generated by $\{\upsilon_1\}$, i.e. $H_1= \C \,
\upsilon_1$, and $K_1:= \upsilon_1^\perp$. Thus $H= H_1 \oplus
K_1$, that is, for each $x \in H$, we write $x= x_1 + y_1$.
Moreover, by Remark \ref{RRT}, we have
$
  T (\upsilon_1^\perp) \subset \upsilon_1^\perp,
$
and setting $\beta_1:= \|T\|$, we could write
$$
  \begin{aligned}
  T x&= \beta_1 \langle x, \upsilon_1 \rangle \, \upsilon_1
  + T y_1
  = \beta_1 \, (\upsilon_1 \otimes \upsilon_1) \, x + T y_1.
  \end{aligned}
$$
Therefore, denoting $T y_1= T_1 x$, it follows that $
  T= \beta_1 \, (\upsilon_1 \otimes \upsilon_1) + T_1,
$ where $T_1$ has the following representation,
$$
  T_1= \left (
  \begin{matrix}
  0  &   0   \\
  0  &   T^1
  \end{matrix}
  \right ),
$$
and $T^1$ is the restriction of $T$ to the subspace $K_1$, i.e.
$T^1= T|_{K_1}$. Now, since $T^1$ is positive and satisfies
$\clg{N}$, hence there exists $\upsilon_2 \in (S \cap K_1)$, such
that $T^1 \upsilon_2= \beta_2 \, \upsilon_2$, where $
  \beta_2= \|T^1\| \leq \beta_1.
$ It is clear that, $\upsilon_1 \perp \upsilon_2$, i.e they are
orthogonal. Analogously, we set $H_2= \C \, \upsilon_2$, that is
the subspace of $H_1$ generated by $\{\upsilon_2\}$, and $K_2=
\upsilon_2^\perp$. Then, $T$ could be written as
$$
  T= \beta_1 \; \upsilon_1 \otimes \upsilon_1 + \beta_2 \; \upsilon_2 \otimes \upsilon_2 + T_2,
$$
with $T_2$ given by
$$
  T_2= \left (
  \begin{matrix}
  0  &   0   &  0   \\
  0  &   0   &  0   \\
  0  &   0   &  T^2
  \end{matrix}
  \right ),
$$
where $T^2$ is the restriction of $T^1$ to the subspace $K_2$,
i.e. $T^2= T^1|_{K_2}$.
Thus, continuing in this way, the operator $T$, which satisfies
$\clg{AN}$, could be written as
\begin{equation}
\label{Decomp}
  T= \sum_{n \geq 1} \beta_n \; \upsilon_n \otimes \upsilon_n + R_1,
\end{equation}
where $\{\upsilon_n\}$ is an orthonormal sequence of  vectors in
$H$, $\{\beta_n\}$ is a decreasing sequence of positive real
numbers, such that, for all $j \geq 1$ $
  T \upsilon_j= \beta_j \; \upsilon_j,
$ and $R_1$ is a remainder operator, satisfying $\|R_1\| \leq
\beta_j$, $(j \geq 1)$.

If $R_1= 0$, then the operator $T$ has the simplest diagonal form.
For instance, it happens when $\dim T(H) < \infty$ or $T$ is
compact.

%
%
\end{pf}

\section*{References}

\end{document}